\documentclass[12pt]{article}
\usepackage{amsmath}
\usepackage{amssymb}
\usepackage{amsmath,amssymb,amsbsy,amsfonts,amsthm,latexsym,
amsopn,amstext,amsxtra,euscript,amscd}

\begin{document}

\def\A{\mathbb{A}}
\def\B{\mathbf{B}}
\def \C{\mathbb{C}}
\def \F{\mathbb{F}}
\def \K{\mathbb{K}}

\def \Z{\mathbb{Z}}
\def \P{\mathbb{P}}
\def \R{\mathbb{R}}
\def \Q{\mathbb{Q}}
\def \N{\mathbb{N}}
\def \Z{\mathbb{Z}}

\def\B{\mathcal B}
\def\e{\varepsilon}

\def\cA{{\mathcal A}}
\def\cB{{\mathcal B}}
\def\cC{{\mathcal C}}
\def\cD{{\mathcal D}}
\def\cE{{\mathcal E}}
\def\cF{{\mathcal F}}
\def\cG{{\mathcal G}}
\def\cH{{\mathcal H}}
\def\cI{{\mathcal I}}
\def\cJ{{\mathcal J}}
\def\cK{{\mathcal K}}
\def\cL{{\mathcal L}}
\def\cM{{\mathcal M}}
\def\cN{{\mathcal N}}
\def\cO{{\mathcal O}}
\def\cP{{\mathcal P}}
\def\cQ{{\mathcal Q}}
\def\cR{{\mathcal R}}
\def\cS{{\mathcal S}}
\def\cT{{\mathcal T}}
\def\cU{{\mathcal U}}
\def\cV{{\mathcal V}}
\def\cW{{\mathcal W}}
\def\cX{{\mathcal X}}
\def\cY{{\mathcal Y}}
\def\cZ{{\mathcal Z}}

\def\f{\frac{|\A||B|}{|G|}}
\def\AB{|\A\cap B|}
\def \Fq{\F_q}
\def \Fqn{\F_{q^n}}

\def\({\left(}
\def\){\right)}
\def\fl#1{\left\lfloor#1\right\rfloor}
\def\rf#1{\left\lceil#1\right\rceil}
\def\Res{{\mathrm{Res}}}

\newcommand{\comm}[1]{\marginpar{
\vskip-\baselineskip \raggedright\footnotesize
\itshape\hrule\smallskip#1\par\smallskip\hrule}}

\newtheorem{lem}{Lemma}
\newtheorem{lemma}[lem]{Lemma}
\newtheorem{cl}{Claim}
\newtheorem{claim}[cl]{Claim}
\newtheorem{prop}{Proposition}
\newtheorem{proposition}[prop]{Proposition }
\newtheorem{thm}{Theorem}
\newtheorem{theorem}[thm]{Theorem}
\newtheorem{cor}{Corollary}
\newtheorem{corollary}[cor]{Corollary}
\newtheorem{prob}{Problem}
\newtheorem{problem}[prob]{Problem}
\newtheorem{ques}{Question}
\newtheorem{question}[ques]{Question}
\newtheorem{rem}{Remark}

\title{On integer values of sum and product of three positive rational numbers}

\author{
{\sc M.~Z.~Garaev} }

\date{}

\maketitle

\begin{abstract}

In 1997 we proved that if $n$ is of the form
$$
4k, \quad 8k-1\quad  {\rm or} \quad 2^{2m+1}(2k-1)+3,
$$
where $k,m\in \mathbb N,$ then there are no positive rational numbers $x,y,z$  satisfying
$$
xyz = 1,  \quad x+y+z = n.
$$
Recently, N. X. Tho proved the following statement: let $a\in\mathbb N$ be odd and let either $n\equiv 0\pmod 4$ or $n\equiv 7\pmod 8$. Then the system of equations
$$
xyz = a,  \quad x+y+z = an.
$$
has no solutions in positive rational numbers $x,y,z.$

A representative example of our result is the following statement: assume that $a,n\in\mathbb N$ are such that at least one of the following conditions hold:

\begin{itemize}
  \item  $n\equiv 0\pmod 4$
  \item   $n\equiv 7\pmod 8 $
  \item   $a\equiv 0\pmod 4$
  \item  $a\equiv 0\pmod 2$ and $n\equiv 3\pmod 4$
  \item  $a^2n^3=2^{2m+1}(2k-1)+27$ for some $k,m\in \mathbb N.$
\end{itemize}
Then the system of equations
$$
xyz = a,  \quad x+y+z = an.
$$
has no solutions in positive rational numbers $x,y,z.$

\end{abstract}

%\paragraph{Mathematical Subject Classification:}

%\paragraph{Keywords:} Diophantine equation, Sylvester transformation, Jacoby symbol.

\section{Introduction}

Let $n\not\in\{0, 3\}$ be an integer. The Diophantine equation
\begin{equation}
\label{eqn:CubeSum}
X^3 + Y^3 + Z^3 = nXYZ,
\end{equation}
has been a subject of much investigation, starting from the works of Sylvester dating back to 1856. Sylvester proved that
if $n=-6,$ then~\eqref{eqn:CubeSum} has no solutions in nonzero integers $X,Y,Z,$ see~\cite{Dic}. The cases $n\in\{-1, 5\}$ were treated
by  Mordell~\cite{Mor}. In 1960 Cassels~\cite{Cas} proved that for $n=1$ the equation~\eqref{eqn:CubeSum} has no nonzero integer solutions.

Further contribution was made by Dofs~\cite{Dof}. One of his results claims that if $n^2+3n+9$ is a prime number and if all prime divisors of $n-3$ are of the form $2\pmod 3$, then~\eqref{eqn:CubeSum} has no solutions in nonzero integers $X,Y,Z$.

Considering the equation
\begin{equation}
\label{eqn:Sier}
\frac{x}{y}+\frac{y}{z}+\frac{z}{x}=n
\end{equation}
for $n=4,$ Sierpinski~\cite[p.80]{Sie} noted that it was not known to him whether this equation has a solutions in positive integers $x,y,z.$ From the mentioned result of Dofs it follows that for $n=4$ the equation~\eqref{eqn:Sier} has no solutions even in nonzero integers $x,y,z$.

In~\cite{Gar1} we proved that if
$$
n\in\{4k, \, 8k-1,\, 2^{2m+1}(2k-1)+3\},
$$
where $k$ and $m$ run through all positive integers, then~\eqref{eqn:CubeSum} has no solutions  in positive integers $X,Y,Z$. As a consequence of the Sylvester transformation, we showed that for the indicated set of values of $n$ there do not exist positive rational numbers $x,y,z$  satisfying
$$
xyz = 1,  \quad x+y+z = n.
$$
In particular, the equation~\eqref{eqn:Sier}  has no positive integer solutions not only for $n=4$ but also for any $n\equiv 0\pmod 4$ and for many odd values of $n$ as well. In~\cite{Gar3} we also proved that
the equation
$$
\Bigr(\frac{x}{y}\Bigl)^u+\Bigr(\frac{y}{z}\Bigl)^v+\Bigr(\frac{z}{x}\Bigl)^w = 4t
$$
has no solutions in positive integers $x,y,z,t,u,v,w.$

Recently, Tho~\cite{Tho} proved the following statement: let $a\in\N$ be odd and let either $n\equiv 0\pmod 4$ or $n\equiv 7\pmod 8$. Then the system of equations
\begin{equation}
\label{eqn:EqTho1}
xyz = a,  \quad x+y+z = an.
\end{equation}
has no solutions in positive rational numbers $x,y,z.$

In the present paper we shall prove the following result.

\begin{theorem}
\label{thm: Main}
Let $a,b$ and $c$ be positive integers such that at least one of the following conditions hold:
\begin{itemize}
\item   $c\equiv 0\pmod 4 $
\item  $bc\equiv 7\pmod 8$
\item  $a\equiv 0\pmod 4$\, and\, $b\equiv 1\pmod 2$
\item  $a\equiv 0\pmod 2$ \,and\, $bc\equiv 3 \pmod 4$
\item   $a\equiv 1\pmod 2,$\, $b\equiv 2\pmod 4$\, and\, $c\equiv 0\pmod 2$
\item  $a^2bc^3 = 2^{2m+1}(2k-1)+27$ \,for some\, $k,m\in \N.$
\end{itemize}
Then there do not exist positive rational numbers $x,y,z$  satisfying
\begin{equation}
\label{eqn:MainMe}
xyz = ab^2,\quad x+y+z = abc.
\end{equation}
\end{theorem}

Taking $b=1$ and $n=c,$ we obtain the following consequence.

\begin{corollary}
\label{cor:Main}
Let  $a$ and $n$ be positive integers such that at least one of the following conditions hold:
\begin{itemize}
  \item  $n\equiv 0\pmod 4$
  \item   $n\equiv 7\pmod 8 $
  \item   $a\equiv 0\pmod 4$
  \item  $a\equiv 0\pmod 2$ and $n\equiv 3\pmod 4$
  \item  $a^2n^3=2^{2m+1}(2k-1)+27$ \,for some\, $k,m\in \N.$
\end{itemize}
Then~\eqref{eqn:EqTho1} has no solutions in positive rational numbers  $x,y,z.$
\end{corollary}

\section{Main Lemma}

The following statement is the key of the proof of Theorem~\ref{thm: Main}.

\begin{lemma}
\label{lem:n2=n}
Let
$$
n\in\{16k-4, \,64k, \,32k-16, \,8k-1, \, 2^{2m+1}(2k-1)+27\},
$$
where $k$ and $m$ run through all positive integers. Then the equation
\begin{equation}
\label{eqn: n2 vs n}
x^3 + y^3 +n^2 z^3 = nxyz
\end{equation}
has no solutions in positive integers $x,y,z.$
\end{lemma}

Lemma~\ref{lem:n2=n} is implicitly contained in our note~\cite{Gar2},  without proof. Here, we shall give a complete proof of this
statement.

We shall often use the following simple observation.

\begin{claim}
\label{cl: pairwise primes}
Let $x,y,z,A,B$ be positive integers such that
$$
x^3+y^3+Az^3 = Bxyz.
$$
Assume that $A$ is cubefree, $(x,y,z)=1,$ and that any prime divisor of $A$ is also a prime divisor of $B.$ Then
$$
(x,y)=(y,z)=(z,x)=(x,A)=(y,A)=1.
$$
\end{claim}
Indeed, if a prime number $p$ divides $(x,y),$ then $p$ also divides $A,$ implying that $p$ divides $B,$ and therefore $p^3$ divides $x^3,y^3, Bxyz.$
Hence $p^3$ divides $Az^3,$ which is impossible as $A$ is cubefree and $(x,y,z)=1.$ Therefore $(x,y)=1,$ from which it also follows the remaining conclusion.

\smallskip

We proceed to prove Lemma~\ref{lem:n2=n}. Assume that~\eqref{eqn: n2 vs n} holds for some positive integers $x,y,z.$ We shall consider five cases.

\subsection{ The case  $n=16k-4$}

\smallskip

We represent $n$ in the form,
$$
n=4P_1P_2^2P_3^3,
$$
where $P_1$ and $P_2$ are relatively prime squarefree positive integers. We also have that
$$
P_1P_2^2P_3^3\equiv 3\pmod 4,
$$
which means that $P_1,P_2,P_3$ are odd and
$$
P_1P_3 \equiv 3 \pmod 4.
$$

We have that
$$
x^3 + y^3 + 16 P_1^2P_2^4P_3^6 z^3 = 4P_1P_2^2P_3^3 xyz.
$$
Therefore,
$$
x^3 + y^3 + 2P_1^2 P_2z_0^3 = 2 P_1P_2P_3 xyz_0,
$$
where $z_0 = 2P_2P_3^2z.$ Letting $d=(x,y,z_0)$ and denoting $x/d, y/d, z_0/d$ again by $x,y,z,$
we get that the equality
$$
x^3 + y^3 + 2P_1^2 P_2z^3 = 2 P_1P_2P_3 xyz
$$
holds for some positive integers $x,y,z$ with $(x,y,z)=1.$ Note that $P_1^2P_2$ is an odd cubefree integer. Thus, from Claim~\ref{cl: pairwise primes}
it follows that $x$ and $y$ are odd integers satisfying
$$
(x,y)=(y,z)=(z,x)=(x, 2P_1P_2)=(y,2P_1P_2)=1.
$$
From
$$
x^3+y^3 = 2P_1P_2z(P_3xy - P_1z^2),
$$
it also follows that
\begin{equation}
\label{eqn: first x+y = 4}
x+y\equiv 0\pmod 4.
\end{equation}
There are two possibilities, depending on whether $z$ is an even or an odd number.

Let $z$ be an even number. Then
$z=2^r z_1,$ where $r\in\N$ and $z_1$ is odd. Then,
$$
x^3\equiv -2^{3r+1}P_1^2P_2z_1^3 \pmod {2^{r+1}P_1P_2P_3xz_1 - y^2},
$$
whence
$$
X^2\equiv -2^{3r+1}P_1^2P_2xz_1^3 \pmod {2^{r+1}P_1P_2P_3xz_1 - y^2}, \quad X=x^2.
$$
Observe that the numbers  $2^{3r+1}P_1^2P_2xz_1^3$ and $2^{r+1}P_1P_2P_3xz_1 - y^2$ are relatively prime positive integers. Indeed, if a prime number $p$ divides both of these numbers, then clearly $p$ divides $y$, and therefore $x$ as well, which contradicts to $(x,y)=1.$
Since $2^{r+1}P_1P_2P_3xz_1 - y^2$ is odd,  we can pass to the Jacoby symbol and get that
\begin{equation}
\label{eqn: first jacob}
1 = \left(\frac{-2^{3r+1}P_1^2P_2xz_1^3}{2^{r+1}P_1P_2P_3xz_1-y^2}\right) = \left(\frac{-2^{3r+1}P_2xz_1}{2^{r+1}P_1P_2P_3xz_1-y^2}\right).
\end{equation}
Note that $y^2\equiv 1\pmod 8,$ which implies that
$$
\left(\frac{-1}{2^{r+1}P_1P_2P_3xz_1-y^2}\right) = -1.
$$
We also have that
$$
\left(\frac{2}{2^{r+1}P_1P_2P_3xz_1-y^2}\right)^{3r+1} = 1.
$$
Indeed, it is trivial for $r=1$, and for $r\ge 2$ it follows from the fact that $2^{r+1}P_1P_2P_3xz_1-y^2\equiv -1\pmod 8.$

Thus, from~\eqref{eqn: first jacob} we get that
\begin{equation}
\label{eqn: first contr jacob}
\left(\frac{P_2xz_1}{2^{r+1}P_1P_2P_3xz_1-y^2}\right) = -1.
\end{equation}
On the other hand, from  $2^{r+1}P_1P_2P_3xz_1-y^2\equiv 3 \pmod 4$ and the quadratic  reciprocity law, we have that
\begin{equation*}
\begin{split}
\left(\frac{P_2xz_1}{2^{r+1}P_1P_2P_3xz_1-y^2}\right) = (-1)^{(P_2xz_1-1)/2}\left(\frac{2^{r+1}P_1P_2P_3xz_1-y^2}{P_2xz_1}\right)\\=
(-1)^{(P_2xz_1-1)/2}\left(\frac{-y^2}{P_2xz_1}\right)= (-1)^{(P_2xz_1-1)/2}\left(\frac{-1}{P_2xz_1}\right)=1,
\end{split}
\end{equation*}
which contradicts to~\eqref{eqn: first contr jacob}.

Therefore, $z$ should be an odd number. Then, from~\eqref{eqn: first x+y = 4} it follows that
either $xz\equiv 3 \pmod 4$ or $yz\equiv 3\pmod 4.$
Without loss of generality, we can suppose that $xz\equiv 3\pmod 4.$ We have that
$$
X^2\equiv -2P_1^2P_2xz^3 \pmod {2P_1P_2P_3xz - y^2}, \quad X=x^2.
$$
Since the numbers $2P_1^2P_2xz^3$ and $2P_1P_2P_3xz - y^2$ are relatively prime positive integers, with $2P_1P_2P_3xz - y^2$ being odd, we can use the Jacoby symbol:

\begin{equation}
\label{eqn: Jac Case 1 z odd}
1= \left(\frac{-2P_1^2P_2xz^3}{2P_1P_2P_3xz-y^2}\right) = \left(\frac{-2P_2xz}{2P_1P_2P_3xz-y^2}\right).
\end{equation}

We recall that $P_1P_3\equiv xz\equiv 3\pmod 4.$ Hence,
$$
2P_1P_2P_3xz-y^2\equiv 2P_2-1\pmod 8.
$$
This and $2P_2-1\equiv 1\pmod 4$ implies that,
$$
\left(\frac{-2}{2P_1P_2P_3xz-y^2}\right)=\left(\frac{-1}{2P_2-1}\right)\left(\frac{2}{2P_2-1}\right)=(-1)^{(P_2^2-P_2)/2}=(-1)^{(1-P_2)/2}.
$$
From $2P_1P_2P_3xz-y^2\equiv 1 \pmod 4,$ $xz\equiv 3\pmod 4$ and the quadratic reciprocity law, we also have that
\begin{equation*}
\begin{split}
\left(\frac{P_2xz}{2P_1P_2P_3xz-y^2}\right) &= \left(\frac{2P_1P_2P_3xz-y^2}{P_2xz}\right)\\&=\left(\frac{-y^2}{P_2xz}\right)=(-1)^{(P_2xz-1)/2}=(-1)^{(3P_2-1)/2}.
\end{split}
\end{equation*}
Putting the last two relations into~\eqref{eqn: Jac Case 1 z odd}, we conclude that
$$
1 = (-1)^{(1-P_2)/2}(-1)^{(3P_2-1)/2}=(-1)^{P_2} = -1.
$$
The obtained contradiction finishes the case $n=16k-4$ of our lemma.

\subsection{ The case  $n=64k$}

We let
$$
n=64P_1P_2^2P_3^3,
$$
where $P_1,P_2$ are squarefree relatively prime integers. Similar to the previous case, we arrive that
there exists positive integers $x,y,z$ such that $(x,y,z)=1$ and
$$
x^3+y^3+P_1^2P_2z^3 = 4P_1P_2P_3xyz.
$$
By Claim~\ref{cl: pairwise primes}, we have that
\begin{equation}
\label{eqn: cl 64k}
(x,y)=(y,z)=(z,x)=(x, P_1P_2)=(y, P_1P_2)=1.
\end{equation}
It follows that at least one of the numbers $x$ and $y$ is an odd number. Assume that $y$ is odd.
Then we can write,
$$
x^3 \equiv -P_1^2P_2z^3 \pmod {4P_1P_2P_3xz-y^2}.
$$
From~\eqref{eqn: cl 64k} it follows that the numbers $P_1^2P_2z^3x$ and $4P_1P_2P_3xz-y^2$ are relatively prime.  Indeed, if a prime number $p$
divides both them, then $p$ divides $y$ and hence, $p$ divides $P_1^2P_2z^3x,$ which contradicts to~\eqref{eqn: cl 64k}.

Thus, the numbers $P_1^2P_2z^3x$ and $4P_1P_2P_3xz-y^2$ are relatively prime positive integers, with $4P_1P_2P_3xz-y^2$ being odd, and therefore we can use the Jacoby symbol in the congruence
$$
X^2 \equiv -P_1^2P_2xz^3 \pmod {4P_1P_2P_3xz-y^2}, \quad X=x^2.
$$
Setting $P_2xz=2^{r}t$, where $t$ is odd and $r$ is a nonnegative integer, we get that
\begin{equation*}
\begin{split}
1&=\left(\frac{-P_1^2P_2xz^3}{4P_1P_2P_3xz-y^2}\right) = \left(\frac{-2^rt}{2^{r+2}P_1P_3t-y^2}\right) \\&= \left(\frac{-1}{2^{r+2}P_1P_3t-y^2}\right)\left(\frac{2}{2^{r+2}P_1P_3t-y^2}\right)^r\left(\frac{t}{2^{r+2}P_1P_3t-y^2}\right)\\&=
-\left(\frac{t}{2^{r+2}P_1P_3t-y^2}\right)=-(-1)^{(t-1)/2}\left(\frac{2^{r+2}P_1P_3t-y^2}{t}\right)\\
&=-(-1)^{(t-1)/2}\left(\frac{-y^2}{t}\right)=-(-1)^{(t-1)/2}\left(\frac{-1}{t}\right)=-1.
\end{split}
\end{equation*}
We have a contradiction, which finishes the case $n=64k.$

\subsection{The case $n =32k-16$}

We write
$$
n=16P_1P_2^2P_3^3,
$$
where $P_1,P_2,P_3$ are odd integers, $P_1$ and $P_2$ are relatively prime squarefree integers. Similar to the previous cases, we
get that for some positive integers $x,y,z$ with $(x,y,z)=1$ one has
$$
x^3+y^3+4P_1^2P_2z^3 = 4P_1P_2P_3xyz.
$$
By Claim~\ref{cl: pairwise primes} we have that
$$
(x,y)=(y,z)=(z,x)=(x, 2P_1P_2)=(y,2P_1P_2)=1.
$$
From our equation it follows that,
$$
x^4 = - 4P_1^2P_2xz^3\pmod {4P_1P_2P_3xz-y^2}.
$$
Since the numbers $4P_1^2P_2xz^3$ and $4P_1P_2P_3xz-y^2$ are relatively prime positive integers, with $4P_1P_2P_3xz-y^2$ being odd, we can use the Jacoby symbol:
$$
1=\left(\frac{-4P_1^2P_2xz^3}{4P_1P_2P_3xz-y^2}\right) = \left(\frac{-P_2xz}{4P_1P_2P_3xz-y^2}\right).
$$
Let
$P_2xz=2^rt$, where $t$ is odd and $r$ is a nonnegative integer. Note that
$$
\left(\frac{-1}{2^{r+2}P_1P_3t-y^2}\right)=-1, \quad \left(\frac{2}{2^{r+2}P_1P_3t-y^2}\right)^r=1.
$$
Hence,
\begin{equation*}
\begin{split}
1&=\left(\frac{-P_2xz}{4P_1P_2P_3xz-y^2}\right)=-\left(\frac{t}{2^{r+2}P_1P_3t-y^2}\right)\\
&=(-1)^{(t-1)/2)}\left(\frac{2^{r+2}P_1P_3t-y^2}{t}\right)
=-(-1)^{(t-1)/2)}\left(\frac{-y^2}{t}\right)\\&=-(-1)^{(t-1)/2)}\left(\frac{-1}{t}\right)=-1.
\end{split}
\end{equation*}

The obtained contradiction finishes the case $n=32k-16.$

\subsection{The case $n = 8k-1$}

We write $n=P_1P_2^2P_3^3,$ where $P_1,P_2,P_3$ are odd positive integers, $P_1$ and $P_2$ are squarefree relatively prime integers. The condition $n\equiv -1\pmod 8$ implies that $P_1P_3\equiv -1\pmod 8.$ As in the previous cases, we get that there are positive integers $x,y,z$
with $(x,y,z)=1$ such that
$$
x^3+y^3+P_1^2P_2z^3 = P_1P_2P_3xyz.
$$
Since $(x,y,z)=1$ and  $P_1^2P_2$ is cubrefree, from Claim~\ref{cl: pairwise primes} we get that
$$
(x,y)=(y,z)=(z,x)=(x, P_1P_2)=(y, P_1P_2)=1.
$$

We distinguish two possibilities, depending on whether $z$ is even or odd. Assume that $z$ is even. Then $xy$ is odd. We also have that
$$
x^4\equiv -xy^3 \pmod {P_1P_3xy-P_1^2z^2}
$$
Since the positive integers $xy^3$ and $P_1P_3xy-P_1^2z^2$ are relatively prime, with $P_1P_3xy-P_1^2z^2$ being odd, we can use the Jacoby symbol. Taking into account that $P_1P_3\equiv -1\pmod 8,$ we consequently get that
\begin{equation*}
\begin{split}
1 &= \left(\frac{-xy^3}{P_1P_3xy-P_1^2z^2}\right) = (-1)^{(P_1P_3xy-1)/2}\left(\frac{xy}{P_1P_3xy-P_1^2z^2}\right)\\
 &=(-1)^{(-xy-1)/2}(-1)^{((xy-1)/2)((-xy-1))/2}\left(\frac{P_1P_3xy-P_1^2z^2}{xy}\right)\\
 &= (-1)^{(-xy-1)/2}\left(\frac{-P_1^2z^2}{xy}\right)=(-1)^{(-xy-1)/2}(-1)^{(xy-1)/2}=-1.
\end{split}
\end{equation*}
Here we used that $(-1)^{((xy-1)/2)((-xy-1))/2}=1.$

The obtained contradiction shows that $z$ should be an odd number. We again distinguish two cases, depending on whether
$xy$ is an even or an odd number.

Assume that $xy$ is an even number. Without loss of generality, we can assume that $y$ is even. Then $x$ is odd. We have that
$$
x^4 \equiv -P_1^2 P_2 x z^3 \pmod {P_1P_2P_3xz-y^2}.
$$
Since $P_1^2 P_2 x z^3$ and $P_1P_2P_3xz-y^2$ are relatively prime odd positive integers, we can use the Jacoby symbol and  apply the quadratic reciprocity law. We get that
\begin{equation}
\label{eqn: 8k-1 u odd}
1= \left(\frac{-P_1^2P_2xz^3}{P_1P_2P_3xz-y^2}\right) =(-1)^{(P_1P_2P_3xz-1)/2}\left(\frac{P_2xz}{P_1P_2P_3xz-y^2}\right).
\end{equation}
Next, we have that
\begin{equation*}
\begin{split}
\left(\frac{P_2xz}{P_1P_2P_3xz-y^2}\right)&=(-1)^{((P_1P_2P_3xz-1)/2) \cdot ((P_2xz-1)/2)}\left(\frac{P_1P_2P_3xz-y^2}{P_2xz}\right)\\&=
(-1)^{((P_1P_2P_3xz-1)/2) \cdot ((P_2xz-1)/2)}\left(\frac{-y^2}{P_2xz}\right)\\&=(-1)^{((P_1P_2P_3xz-1)/2) \cdot ((P_2xz-1)/2)}(-1)^{(P_2xz-1)/2}.
\end{split}
\end{equation*}
Inserting this into~\eqref{eqn: 8k-1 u odd}, we see that the number
$$
u= \frac{P_1P_2P_3xz-1}{2}+\frac{P_1P_2P_3xz-1}{2} \cdot \frac{P_2xz-1}{2}+\frac{P_2xz-1}{2}
$$
should be even. However,  $P_1P_3\equiv -1\pmod 8,$ and therefore,
\begin{eqnarray*}
u\equiv \frac{-P_2xz-1}{2}+\frac{-P_2xz-1}{2} \cdot \frac{P_2xz-1}{2} + \frac{P_2xz-1}{2}\\ \equiv -1-\frac{(P_2xz)^2-1}{4}\equiv 1\pmod 2.
\end{eqnarray*}
Contradiction.

Therefore, we remained with the case when $x,y,z$ are all odd numbers.
From our equation we have that
$$
x^4\equiv -xy^3\pmod {P_2z},\,\, x^{4}\equiv -P_1^2P_2z^3x\pmod y, \,\, y^{4}\equiv -P_1^2P_2z^3y\pmod x.
$$
Hence, using the Jacoby symbol, we get that
$$
\left(\frac{-xy^3}{P_2z}\right)= \left(\frac{-P_1^2P_2z^3x}{y}\right) = \left(\frac{-P_1^2P_2z^3y}{x}\right)=1,
$$
whence,
$$
(-1)^{(P_2z-1)/2}\left(\frac{xy}{P_2z}\right) = (-1)^{(y-1)/2}\left(\frac{P_2zx}{y}\right) = (-1)^{(x-1)/2}\left(\frac{P_2zy}{x}\right)=1.
$$
Taking the product, we get that
\begin{equation}
\label{eqn: 8k-1 xyz odd}
(-1)^{(P_2z-1)/2 + (x-1)/2 +(y-1)/2}\left(\frac{xy}{P_2z}\right)\left(\frac{P_2zx}{y}\right)\left(\frac{P_2zy}{x}\right) = 1.
\end{equation}
Furthermore, from the properties of the Jacoby symbol and the quadratic reciprocity law, we have that
\begin{equation*}
\begin{split}
\left(\frac{xy}{P_2z}\right)&\left(\frac{P_2zx}{y}\right)\left(\frac{P_2zy}{x}\right)=
\left(\frac{x}{P_2z}\right)\left(\frac{y}{P_2z}\right)\left(\frac{P_2z}{y}\right)\left(\frac{x}{y}\right)\left(\frac{P_2z}{x}\right)\left(\frac{y}{x}\right)
\\&= \left(\frac{x}{P_2z}\right)\left(\frac{P_2z}{x}\right)\left(\frac{y}{P_2z}\right)\left(\frac{P_2z}{y}\right)\left(\frac{x}{y}\right)\left(\frac{y}{x}\right)
\\&=(-1)^{((x-1)/2)((P_2z-1)/2)+((P_2z-1)/2)((y-1)/2)+((y-1)/2)((x-1)/2)}.
\end{split}
\end{equation*}
Hence, inserting this into~\eqref{eqn: 8k-1 xyz odd}, we obtain that the number
$$
\frac{x-1}{2}+\frac{y-1}{2}+\frac{P_2z-1}{2}+\frac{x-1}{2}\cdot \frac{y-1}{2}+\frac{y-1}{2}\cdot\frac{P_2z-1}{2}+\frac{P_2z-1}{2}\cdot \frac{x-1}{2}
$$
is an even integer. Note that for integers $u,v,w$ the number
$$
u+v+w+uv+vw+wu
$$ is even if and only if either $u,v,w$ are all even, or are all odd numbers.
Hence, we have
$$
x\equiv y\equiv P_2z\equiv r_0\pmod 4,
$$
where $r_0\in \{1,3\}.$ Take
$$
x= r_0 + 4x_0,\quad y=r_0 + 4y_0, \quad P_2z= r_0 +4z_0.
$$
Returning to our equation and recalling that $P_1P_3\equiv -1\pmod 8$, we see that
$$
x+y+P_2z+P_2xyz\equiv 0\pmod 8.
$$
Thus,
$$
3r_0 +4(x_0+y_0+z_0)+ r_0^3+4r_0^2(x_0+y_0+z_0)\equiv 0\pmod 8.
$$
But this is not true, as the left hand side is
$$
r_0(3+r_0^2)+4(1+r_0^2)(x_0+y_0+z_0)\equiv 4\pmod 8.
$$
This contradiction finishes the case $n=8k-1.$

\subsection{The case $n = 2^{2m+1}(2k-1)+27$}

We have that
$$
n-27\equiv 2^{2m+1}\pmod {2^{2m+2}}.
$$
We  write $n=P_1P_2^2P_3^3,$ where $P_1,P_2,P_3$ are odd positive integers, $P_1,P_2$ are squarefree relatively prime integers.
As in the previous cases, we get that there are positive integers $x,y,z$
with $(x,y,z)=1$ such that
$$
x^3+y^3+P_1^2P_2z^3 = P_1P_2P_3xyz.
$$
Since $(x,y,z)=1$ and  $P_1^2P_2$ is cubefree, from Claim~\ref{cl: pairwise primes} we get that
$$
(x,y)=(y,z)=(z,x)=(x, P_1P_2)=(y,P_1P_2)=1.
$$

Clearly, there exists odd integers $u$ and $v$ such that
$$
P_1\equiv u^3\pmod {2^{2m+2}},\quad P_2\equiv v^3\pmod {2^{2m+2}}.
$$
For instance, one can take $u=P_1^{(2^{2m+1}+1)/3}$ and similarly $v.$
Then,
$$
x^3 +y^3 + u^6v^3z^3 -3xyu^2vz\equiv u^3v^3P_3xyz - 3xyu^2vz \pmod {2^{2m+2}}.
$$
Using the decomposition
$$
4(A^3+B^3+C^3-3ABC)=(A+B+C)\Bigl((2A-B-C)^2+3(B-C)^2\Bigr),
$$
we get that
\begin{equation*}
\begin{split}
(x+y+u^2vz)\Big((&2x-y-u^2vz)^2 +3(y-u^2vz)^2\Big)\\ &\equiv 4xyzu^2v(uv^2P_3-3)\pmod {2^{2m+4}}.
\end{split}
\end{equation*}
Multiplying  both hand side by $(uv^2P_3)^2+3(uv^2P_3)+9$ we obtain that
\begin{equation*}
\begin{split}
(x+y+u^2vz)\Big((2x-y-u^2vz)^2 &+3(y-u^2vz)^2\Big)\Big((uv^2P_3)^2+3(uv^2P_3)+9\Big)\\ \equiv 4xyzu^2v(u^3v^6P_3^3-&27)\equiv 4xyzu^2v(P_1P_2^2P_3^3-27)\\
\equiv 4xyzu^2v&(n-27)\equiv 2^{2m+3}xyz\pmod {2^{2m+4}}.
\end{split}
\end{equation*}
Since $(uv^2P_3)^2+3(uv^2P_3)+9$ is odd, we have that
$$
(x+y+u^2vz)\Big((2x-y-u^2vz)^2 +3(y-u^2vz)^2\Big)\equiv 2^{2m+3}xyz\pmod {2^{2m+4}}.
$$
If $xyz$ were odd, we would have that $x+y+u^2vz$ is odd, implying that
$$
(2x-y-u^2vz)^2 +3(y-u^2vz)^2\equiv 2^{2m+3}\pmod {2^{2m+4}}.
$$
This is obviously impossible, as the highest power of $2$ which divides a number of the form $r^2+3s^2$ should be even.

Thus, we proved that $xyz$ is an even number.
Since $x,y,z$ are pairwise primes, it follows that exactly one of them is even.

Now we show that $z$ should be even. Assume contrary, let $z$ be odd. Then, exactly one of the numbers $x$ and $y$ is even. Without loss of generality, we can assume that $y$ is even, $x$ is odd. From our equation we have the congruence
$$
x^4\equiv -P_1^2P_2z^3x \pmod {P_1P_2P_3xz-y^2}.
$$
Since $P_1^2P_2z^3x$ and $P_1P_2P_3xz-y^2$ are relatively prime positive integers, and $P_1P_2P_3xz-y^2$ is odd, we can use the Jacoby symbol.
We deduce that
\begin{equation*}
\begin{split}
1 &= \left(\frac{-P_1^2P_2z^3x}{P_1P_2P_3xz-y^2}\right)=(-1)^{(P_1P_2P_3xz-1)/2}\left(\frac{P_2zx}{P_1P_2P_3xz-y^2}\right)\\
&=(-1)^{(P_1P_2P_3xz-1)/2} (-1)^{((P_2zx-1)/2)((P_1P_2P_3xz-1)/2)}\left(\frac{P_1P_2P_3xz-y^2}{P_2zx}\right)\\&=
(-1)^{(P_1P_2P_3xz-1)/2} (-1)^{((P_2zx-1)/2)((P_1P_2P_3xz-1)/2)}\left(\frac{-y^2}{P_2zx}\right)\\&=
(-1)^{(P_1P_2P_3xz-1)/2} (-1)^{((P_2zx-1)/2)((P_1P_2P_3xz-1)/2)}(-1)^{((P_2zx-1)/2)}\\&=
(-1)^{r+s+rs},
\end{split}
\end{equation*}
where
$$
r=\frac{P_1P_2P_3xz-1}{2}, \quad s = \frac{P_2zx-1}{2}.
$$
Thus, $r+s+rs$ should be an even number, implying that both $r$ and $s$ are even numbers. That is,
$$
P_1P_2P_3xz\equiv 1\pmod 4,\quad P_2xz\equiv 1\pmod 4.
$$
It follows that $P_1P_3\equiv 1\pmod 4,$ whence we get that
$$
n=P_1P_2^2P_3^3\equiv P_1P_3 \equiv 1\pmod 4,
$$
which is contradiction with $n\equiv 3\pmod 4.$

Thus, we get that $z$ is an even number, $xy$ is odd. From our equation we have that
$$
x^4\equiv -xy^3 \pmod {P_1P_3xy-P_1^2z^2}.
$$
Again, $xy^3$ and  $P_1P_3xy-P_1^2z^2$ are relatively prime positive integers, with latter being odd.
Hence, we can use the Jacoby symbol and get that
\begin{equation*}
\begin{split}
1 &= \left(\frac{-xy^3}{P_1P_3xy-P_1^2z^2}\right)=(-1)^{(P_1P_3xy-1)/2}\left(\frac{xy}{P_1P_3xy-P_1^2z^2}\right)\\
&=(-1)^{(P_1P_3xy-1)/2} (-1)^{((xy-1)/2)((P_1P_3xy-1)/2)}\left(\frac{P_1P_3xy-P_1^2z^2}{xy}\right)\\&=
(-1)^{(P_1P_3xy-1)/2} (-1)^{((xy-1)/2)((P_1P_3xy-1)/2)}\left(\frac{-P_1^2z^2}{xy}\right)\\&=
(-1)^{(P_1P_3xy-1)/2} (-1)^{((xy-1)/2)((P_1P_3xy-1)/2)}(-1)^{(xy-1)/2}\\&=
(-1)^{r+s+rs},
\end{split}
\end{equation*}
where
$$
r=\frac{xy-1}{2}, \quad s=\frac{P_1P_3xy-1}{2}.
$$
Thus, $r+s+rs$ should be an even number, implying that $r$ and $s$ are both even numbers. It follows that
$$
xy\equiv 1\pmod 4,\quad P_1P_3xy\equiv 1\pmod 4,
$$
implying again that $P_1P_3 1\pmod 4$ and
$$
n=P_1P_2^2P_3^3\equiv P_1P_3 \equiv 1\pmod 4.
$$
This contradiction finishes the case $n=2^{2m+1}(2k-1)+27.$

Lemma~\ref{lem:n2=n} is proved.

\section{Proof of the theorem}

We shall use the following statement due to Sylvester~\cite{Dic}.

\begin{lemma}
\label{lem: sylv transf}
Let $A, B, C, D, \alpha, \beta, \gamma$ be arbitrary real numbers for which
$$
A\alpha^3+B\beta^3+C\gamma^3 = D\alpha\beta\gamma.
$$
Then
$$
f^3 + g^3 + ABCh^3 =Dfgh,
$$
where
\begin{equation*}
\begin{split}
f&=A^2B\alpha^6\beta^3 + B^2C\beta^6\gamma^3+C^2A\gamma^6\alpha^3-3ABC\alpha^3\beta^3\gamma^3,\\
g&=AB^2\alpha^3\beta^6 + BC^2\beta^3\gamma^6+CA^2\gamma^3\alpha^6-3ABC\alpha^3\beta^3\gamma^3,\\
h&=\alpha\beta\gamma\left(A^2\alpha^6+B^2\beta^6+C^2\gamma^6-AB\alpha^3\beta^3-BC\beta^3\gamma^3-CA\gamma^3\alpha^3\right).
\end{split}
\end{equation*}
\end{lemma}

We proceed to prove our theorem. Assume that $a,b,c$ satisfies one of the conditions listed in the theorem and assume that the positive rational numbers $x,y,z$
are such that
$$
xyz=ab^2,\quad x+y+z = abc
$$
First of all we note that if $x=y=z,$ then $a^2bc^3=27$ which is not satisfied by conditions of the theorem.

We apply Lemma~\ref{lem: sylv transf} with
$$
A=x,\, B=y,\, C=z, \quad \alpha=\beta=\gamma=1, \quad D =abc.
$$
It follows that
$$
f^3+g^3+ab^2 h^3 = abcfgh,
$$
where
\begin{equation*}
\begin{split}
f=x^2y+y^2z&+z^2x-3xyz,\quad g = xy^2+yz^2+zx^2-3xyz, \\
 &h = x^2+y^2+z^2-xy-yz-zx.
\end{split}
\end{equation*}
Since not all the positive rational numbers $x,y,z$ are equal, from the elementary inequalities we get that $f,g$ and $h$ are also positive rational numbers.
We have that
$$
f^3+g^3+(a^2bc^3)^2h_1^3 = a^2bc^3 fgh_1,
$$
where $h_1= h/ac^2.$ Thus, for $n=a^2bc^3$ the equation
$$
x^3+y^3+n^2z^3 = n xyz,
$$
has a solution in positive rational numbers $x,y,z$. Therefore, it also has a solution in positive integers $x,y,z.$

If $c\equiv 0\pmod 4,$ then $n$ is of the form $64k.$

If $bc\equiv 7\pmod 8,$ then $bc^3\equiv 7\pmod 8.$ Hence, $a^2bc^3$ is one of the forms $64k,$ $32k-16,$ $16k-4$ or $8k-1.$

If $a\equiv 0\pmod 4$ and $b\equiv 1\pmod 2,$ then $n$ is either of the form $64k$ or of the form $32k-16.$

If $a\equiv 0\pmod 2$ and $bc\equiv 3\pmod 4,$ then $n$ is one of the forms $64k,$  $32k-16$ or $16k-4.$

If $a\equiv 1\pmod 2,$ $b\equiv 2\pmod 4$ and $c\equiv 0\pmod 2$, then $n$ is either of the form $64 k$ or of the form $32k-16.$

By our Lemma~\ref{lem:n2=n} non of theses cases, neither the case $a^2bc^3=2^{2m+1}(2k-1)+27$ is possible.

Theorem~\ref{thm: Main} is proved.

\section{Remarks}

The proof of  Lemma~\ref{lem:n2=n} also leads to a complete proof of the result that we announced in our work~\cite{Gar2}. That is, if
$$
n\in\{16k-4, \,64k, \,32k-16, \,8k-1, \, 2^{2m+1}(2k-1)+27\},
$$
where $k$ and $m$ run through all positive integers, then the equation
\begin{equation}
\label{eqn: Guy}
\frac{(x+y+z)^3}{xyz}=n
\end{equation}
has no solutions in positive integers $x,y,z.$ Indeed, assume that there are positive integers $x,y,z$ satisfying this equality. We apply the Sylvester transformation to the equality
$$
x\cdot 1^3+y\cdot 1^3+z\cdot 1^3 =(x+y+z)\cdot 1\cdot 1\cdot 1.
$$
It follows that
$$
f^3+g^3+xyzh^3 =(x+y+z)fgh,
$$
for some positive integers $f,g,h.$ Multiplying  by $n^3$ and denoting
$$
X=nf, \quad Y=ng, \quad Z= (x+y+z)h
$$
we get
$$
X^3 +Y^3 + n^2Z^3 = n XYZ.
$$
By Lemma~\ref{lem:n2=n} this is impossible. Thus, for the indicated set of values of $n$ the equation~\eqref{eqn: Guy} has no solutions
in positive integers $x,y,z.$

The problem of representability of integers $n$ in the form~\eqref{eqn: Guy} with some positive integers $x,y,z$ is due to Richard Guy,
see for the details the work of Tho~\cite{Tho1} and the references therein.

Address of the author:\\

M.~Z.~Garaev, Centro de Ciencias Matem\'{a}ticas,  Universidad
Nacional Aut\'onoma de M\'{e}xico, C.P. 58089, Morelia,
Michoac\'{a}n, M\'{e}xico,

Email: {\tt garaev@matmor.unam.mx}

\end{document}